\documentclass{article}

\usepackage[latin1]{inputenc}
\usepackage[english]{babel}
\usepackage{amsfonts}
\usepackage{amssymb}
\usepackage{amsmath}
\usepackage{xypic}
\usepackage[all]{xy}
\usepackage{mathrsfs}
\usepackage{latexsym}
\usepackage{amscd}
\usepackage{theorem}
\usepackage{QEDeng}

\usepackage[mathscr]{eucal}

\usepackage{setspace}                    

\usepackage{tocbibind}

\usepackage{hyperref}

\newtheorem{prop}{Proposition}[section]
\newtheorem{lem}[prop]{Lemma}
\newtheorem{cor}[prop]{Corollary}
\newtheorem{thm}[prop]{Theorem}
\newtheorem{dfn}[prop]{Definition}

\theorembodyfont{\rm}
\newtheorem{rmk}[prop]{Remark}

\newcommand{\qq}{\mathbb{Q}}

\newcommand{\cc}{\mathbb{C}}
\newcommand{\kk}{\Bbbk}

\newcommand{\cou}{{C}}

\newcommand{\gr}{\mathbb{G}\mathrm{r}}

\newcommand{\oo}{\mathcal{O}}
\newcommand{\fais}{\mathcal{F}}

\newcommand{\tais}{\mathcal{T}}

\newcommand{\ii}{\mathcal{I}}

\newcommand{\deriv}[1]{\mathcal{D}(#1)}
\newcommand{\dboun}[1]{{\mathcal{D}}^b(#1)}

\newcommand{\chow}{\mathrm{CH}}
\newcommand{\chowqx}{\mathrm{CH}^*_{\qq}(X)}

         \newif\ifpdf
        \ifx\pdfoutput\undefined
        \pdffalse 
        \else
        \pdfoutput=1 
        \pdftrue
        \fi

        \ifpdf
        \usepackage[pdftex]{graphicx}
        \else
        \usepackage{graphicx}
        \fi

\newcommand{\grafica}{
        \ifpdf
        \DeclareGraphicsExtensions{.pdf, .jpg}
        \else
        \DeclareGraphicsExtensions{.eps, .jpg}
        \fi}
\usepackage{mathpazo}

\title{Stable vector bundles as generators of the  Chow ring}
\author{Ernesto C. Mistretta}
\date{}

\begin{document}

\grafica

%


\maketitle

\begin{abstract}

In this paper we show that the family of stable vector bundles 
gives a set of generators for the Chow ring, 
the $K$-theory and the derived category 
of any smooth projective variety.

\end{abstract}

\section{Introduction}

Let $X$ be a smooth projective variety 
over an algebraically closed field $\Bbbk$,
with a fixed polarization $H$.

The first main result of this note shows that the ideal sheaf $\ii_Z$ 
of an effective cycle $Z \subset X$ 
admits a resolution by polystable vector bundles.
In particular, this shows that the rational Chow ring
$\chowqx$, the $K$-theory $K(X)$, 
and the derived category $\deriv{X}$ 
are generated (in a sense that we will specify) by stable vector bundles.

Note  that it is easy to see that Chern classes of stable 
not necessary locally free sheaves generate 
$\chowqx$ or $K(X)$ (cf. Remark \ref{notfree}).
Since polystability for vector bundles on complex varieties 
is equivalent to the existence of  Hermite-Einstein metrics, 
it seems desirable to work with the more restrictive class of 
locally free stable sheaves.

In the case of a $K3$-surface our result can be compared with a recent
article of Beauville and Voisin. 
In \cite{b-v} they show that all points lying on any rational curve 
are rationally equivalent, 
hence giving rise to the same class $c_X \in \chow (X)$, 
and that $c_2 (X)$ and the intersection product of two 
Picard divisors are multiples of that class.

As the tangent bundle $T_X$ and line bundles are stable, 
one might wonder what happens if we allow arbitrary stable bundles.
Our result shows that second Chern classes of stable bundles 
generate (as a group) 
$\chow^2(X)$, and that this is true on every surface.

Related results, 
using the relation between moduli spaces and Hilbert
schemes (cf. \cite{g-h}), 
and between Hilbert schemes and the second Chow group
(cf. \cite{catr}), 
had been obtained before.

We will first show the main theorem in the case of surfaces, 
as it already gives the above description for  $\chow^2(X)$. 
The higher dimensional case is a generalization of this argument.

\subsection{Notations }

Let $X$ be a smooth projective variety of dimension $n$ 
over an algebraically closed field $\kk$.
The  Chow ring  $\chow^*(X) = \bigoplus C^p(X) /  \sim$ 
is the group of cycles modulo rational equivalence 
graded by codimension. 
Using intersection product of cycles it becomes a commutative graded ring.

For any vector bundle (or a coherent sheaf) $\fais$, Chern classes  
\hbox{$c_i (\fais) \in \chow^i(X)$} and $c(\fais) = \sum c_i(\fais)$ d
efine elements in 
$\chow^* (X)$.

If  $Z$ is an effective cycle, 
we will identify it by abuse of notation with any closed subscheme of 
$X$ having $Z$ as support,
and denote $\ii_Z$ its ideal sheaf.
If $Z$  is an hypersurface 
 $\ii_Z $ is an invertible sheaf and, in particular, stable.
Therefore we will consider only the case 
where $\mathrm{codim} _X Z \geq 2$.

In this paper, 
stability will always mean slope stability 
with respect to the fixed polarization $H$.
Since stability with respect to $H$ or to a multiple $mH$ are
equivalent, 
we can suppose that 
$H$ is sufficiently positive.

\subsection{Acknowledgments}
I wish to thank my advisor Daniel Huybrechts
for many helpful and pleasant discussions,
Claire Voisin  her comments,
and my colleagues of Paris and Rome for their friendly advices.

\section{Zero dimensional cycles on a surface}
\label{surface}

Throughout this section $X$ will be a smooth projective surface,  
$Z$ a 0-dimensional 
subscheme of $X$,
and $\cou \in |H|$  a fixed smooth curve such that  
$ \cou  \cap Z = \emptyset $.
As $H$ is very positive, we suppose $g(\cou) \geqslant 1$.

We will show the following
\begin{prop}
\label{surf}
If $m \gg 0$, and if
\hbox{$V \subset H^0 ( X , \ii_Z (mH))$}
 is a generic subspace of dimension 
$h^0(\cou , \oo_{\cou} (mH))$,
then the sequnce
\begin{equation}
0 \to \ker (ev) \to V \otimes \oo_X \to^{\hspace{-0.35cm}ev}  
~\ii_Z (mH) \to 0 
\end{equation}
is exact and defines a stable vector bundle $M_{Z,m} := \ker (ev)$.
\end{prop}

\subsection{Proof of the proposition}

We remark that if a subspace  $V \subset H^0 ( X , \ii_Z (mH))$
generates $\ii_Z (mH)$, the exact sequence 
\begin{equation}
\label{noy}
 0 \to M \to V \otimes \oo_X \to \ii_Z (mH) \to 0 
\end{equation}
defines a vector bundle $M$, for $X$ has cohomological dimension $2$.

We remark that a sheaf on an arbitrary projective variety is stable 
if its restriction to a hypersurface linearly equivalent to $H$ is stable.

So it is sufficient to show that the restriction 
of $M$ to the curve $\cou$ is a stable vector bundle.
As the chosen curve $\cou$ doesn't intersect  $Z$, 
the restriction of (\ref{noy}) to 
$\cou$ yields a short exact sequence:
\begin{equation}
\label{noc}
 0 \to M |_{\cou} \to V \otimes \oo_{\cou} \to \oo_{\cou} (mH) \to 0  ~. 
\end{equation}

We want to choose the space $V$ so that the sequence  (\ref{noc}) equals:
\begin{equation}
\label{nos}
0 \to M_{\oo_{\cou} (mH)} \to H^0(\cou, \oo_{\cou} (mH)) \otimes \oo_{\cou} \to \oo_{\cou} (mH) \to0 ~. 
\end{equation}

In this case, by general results  (cf. \cite{but}, and Theorem \ref{but} in this paper), the vector bundle 
\hbox{$M |_{\cou} = M_{\oo_{\cou} (mH)}$}
is stable for $m \gg 0$.

We will use the following lemmas:
\begin{lem}
\label{restr}
For  $m \gg 0$, the restriction map 
$$H^0( X , \ii_Z (mH)) \to  H^0(\cou , \oo_{\cou} (mH))$$
induces an isomorphism between a generic subspace  
$V \subset  H^0 ( X , \ii_Z (mH))$ 
of dimension $ h^0(\cou , \oo_{\cou} (mH))$,  
and $H^0(\cou , \oo_{\cou} (mH))$.
\end{lem}

\begin{sloppypar}
\begin{proof}
This follows immediately from the vanishing of 
\hbox{$H^1 (X , \ii_Z ((m-1)H))$} for 
\hbox{$m \gg 0$,} 
and from the consideration that, in the grassmanian 
$Gr ( h^0(\cou , \oo_{\cou} (mH)), H^0 (X , \ii_Z (mH)) )$, 
the spaces $V$ avoiding the subspace $H^0 (X, \ii_Z ((m-1)H))$ 
form an open subset, and project  isomorphically on  
$H^0 (\cou,  \oo_{\cou} (mH))$.

\end{proof}
\end{sloppypar}

So if we show that such a space generates  $ \ii_Z (mH)$, 
then the sequence (\ref{noy}) restricted 
to the curve will give the sequence (\ref{nos}).

Since the dimension $h^0 ( \cou, \oo_{\cou }(mH))$ of such $V$ grows
linearly in $m$,   
this is a consequence of a general lemma which is true 
for a projective variety of any dimension:

\begin{lem}
\label{generate}

Let $Y$ be a projective variety of dimension $n$, 
$E$ a vector bundle of rank  $r$ globally generated on $Y$, 
$\fais$ a coherent sheaf on $Y$, and $H$ an ample divisor. Then:

\begin{enumerate}

\item  If $W \subset H^0 (Y , E)$ 
is a generic subspace of  dimension at least $r+n$, 
then $W$ generates $E$;

\item There are two integers $R, m_0 \geqslant 0$, 
depending on $Y$ and $\fais$, such that for any $m \geqslant m_0$, 
if \hbox{ $V \subset H^0 (Y , \fais (mH))$ }
 is a generic subspace of dimension at least $R$, 
then $V$ generates $\fais (mH)$.

\end{enumerate}

\end{lem}

\begin{proof}

\emph{i.}  Let $W \subset H^0(Y, E)$ be a generic subspace of dimension $v$.
Then the closed subscheme $Y_s \subset Y$ where the evaluation
homomorphism 
$W \otimes \oo_Y \to E$ has rank less than or equal to $s$ is either
empty or of codimension $(v-s)(r-s)$ 
 (cf. \cite{h-l} ch.5, p.121)\footnote{
in \cite{h-l} is used a transversal  version of Kleiman's theorem
which works only in caracteristic 0, 
but the dimension count we need is true in any
characteristic (see \cite{klei}).  
}.
Hence, taking $v= \dim W \geqslant r+n$, and $s = r-1$, 
we see that the evaluation map
must be surjective.

\vspace{.25cm}

\emph{ii.}  By Serre's theorem  there exists a $m_1 \geqslant 0$ 
such that $\fais (m H)$ is globally
generated and acyclic for any $m \geqslant m_1$. 
Hence, there exists a (trivial) globally generated vector bundle $E$ of rank 
$r = h^0 (Y , \fais (m_1 H))$
and a surjection 
\hbox{$ E \twoheadrightarrow \fais (m_1 H) $;}
if we call $\mathcal{K}$ its kernel, 
then $\mathcal{K} (mH)$ is
globally generated and acyclic for any $m \geqslant m_2$, 
and  we have for all $m \geq m_1 + m_2$ :
$$ 0 \to H^0(Y, \mathcal{K} ( (m - m_1)H)) 
\to H^0(Y,  E ((m - m_1)H))
\to H^0 (Y,\fais (mH)) \to 0 ~ .$$

\begin{sloppypar}
Let now  $\nu $ be an integer such that $r+n \leq \nu \leq h^0(Y,\fais (mH))$.
In $Gr( \nu , H^0 (Y, E((m - m_0)H)))$ there is the open subset 
of the spaces $W$ avoiding
$H^0(Y, \mathcal{K} ( (m - m_0)H))$, 
and this open set surjects to $Gr( \nu , H^0 (Y, \fais (mH)))$.

So  a generic $V \subset H^0 (Y, \fais (mH))$ 
of dimension $\nu$ lifts to a generic
\hbox{$W \subset H^0 (Y, E((m - m_0)H))$ }
of dimension $\nu$,  and since $\nu \geq r+n$, 
the first part of this lemma gives the result.
\end{sloppypar}
\end{proof}

Lemmas \ref{restr} and \ref{generate} immediately yield Proposition \ref{surf}.

\subsection{The Chow group of a surface}

We have shown that any effective $0$-cycle $Z$ admits a resolution 
\begin{equation}
\label{noyzm}
 0 \to M_{Z,m} \to V \otimes \oo_X \to \ii_Z(mH) \to 0~,
\end{equation}
where $M_{Z,m}$ is stable and locally free.

\begin{cor}
\label{chowgr}
The Chow group $\chow^2(X)$ is generated as a group by 
$$\{c_2(M) ~|~ \textrm{M is a stable vector bundle} \} ~.$$
\end{cor}

\vspace{.1cm}

\begin{proof}
The class of $Z$ in $\chow^2(X)$ is given by
$ [Z] = -c_2 ( \oo_Z)$,
 hence, 
$$c_2(\ii_Z) = [Z] ~ ,$$ 
furthermore we know that $c_1 (\oo_Z) = c_1 (\ii_Z) = 0$.

Using the sequences
$$ 0 \to \oo_X (-H) \to \oo_X \to \oo_{\cou} \to 0 ~~ \textrm{and,} $$
$$0 \to \ii_Z ((m-1)H) \to \ii_Z (mH) \to \oo_{\cou} (mH) \to 0 ~ , $$
 we can easily calculate the Chern classes  appearing in (\ref{noyzm}):
$$c_1 ( \ii_Z (mH)) = mH ~~ \textrm{and} ~~
c_2 ( \ii_Z (mH)) = c_2 ( \ii_Z ) = [Z] ~.$$

So by the sequence (\ref{noyzm}) we obtain
$$c_1 (M_{Z,m}) = - c_1 (\ii_Z (mH)) = - mH ~ ,~\textrm{and}$$ 
$$ [Z] = c_2 (\ii_Z (mH)) = -c_2 (M_{Z,m}) + m^2 H^2  ~ ,$$
thus second Chern classes of stable vector bundles 
and the class of $H^2$ generate the second Chow group of the surface. 

Clearly, $H^2 = c_2 ( H \oplus H )$ is the second Chern class 
of a polystable vector bundle, 
but it can also be obtained as a linear combination of 
$c_2 (E_i)$ with $E_i$ stable: 
since $H^2 = [Z^{\prime}]$ is an effective cycle, 
we deduce from (\ref{noy}) that
$$[Z^{\prime}] = -c_2 (M_{ Z^{\prime}, m}) + m^2 [Z^{\prime}]  
~~ \textrm{or, equivalently}$$
$$(m^2 -1) H^2 = (m^2 - 1) [Z^{\prime}] = c_2 (M_{Z^{\prime},m})$$
for every $m \gg 0$. 
Choosing $m_1$ and $m_2$ such that $(m_1^2 - 1)$ and $(m_2^2 - 1)$ 
are relatively prime, 
we find that $H^2$ is contained in the subgroup 
of $\chow^2 (X)$ generated by second Chern classes of stable vector bundles.

\end{proof}

\begin{rmk}

This result can also be proven (when $ \mathrm{char} (\kk)=0$)
by using the fact that, 
for every $r>0$, $c_1$, and $c_2 \gg 0$,
stable locally free sheaves form an open dense subset $U$ 
in  the moduli space 
$N = N ( r, c_1, c_2)$
 of semi-stable not necessarily locally free sheaves 
with fixed rank and homological Chern classes
(see \cite{og}).

For any such $N$,
up to desingularizing compactifying and passing to a finite
covering, 
we obtain a homomorphism 
$\phi_{c_2}: CH_0 (N) \to CH_0 (X) $,
which  associates
the class of a point $E \in N$ to the class $c_2 (E) \in CH_0 (X)$.
This morphism is given by the correspondance $c_2(F)$, where $F$ is the 
universal sheaf on $N \times X$.

Next we notice that 
$CH_0 (U)$ spans $CH_0 (N)$: 
in fact if we consider a point
$x \in N$, we can take a curve passing through $x$ and $U$. In the
normalization of this curve, we see that the class of $x$ 
is the
difference of two very ample divisors, so $x$ is rationally
equivalent to a $0$-cycle supported on $X \cap U$. 

Hence the image of the map 
$\phi_{c_2}: CH_0 (N) \to CH_0 (X) $
is spanned by the image of $CH_0 (U)$, and letting vary 
$r, c_1$, and $c_2  \gg 0$, we get the result.

(This remark is due to Claire Voisin).

\end{rmk}

%


\subsection{Bounded families of stable vector bundles generating 
the Chow group of of a surface}

Corollary \ref{chowgr} is interesting in the case of a K3 surface over $\cc$, 
where $\chow(S) = \mathbb{Z} \oplus \mathrm{Pic}(S) \oplus \chow^2 (S)$,
$\mathrm{Pic}(S)$ is a lattice,
and $\chow^2(S)$ is very big (cf. \cite{mum}) and torsion free 
(since $\chow^2(S)_{\textsl{tor}} \subset
\textrm{Alb}(S)_{\textsl{tor}}$ 
for \cite{roit},
and $\textrm{Alb}(S) = 0$). 

Beauville and Voisin have shown in \cite{b-v} that every point lying
on a rational curve has the same class $c_S \in \chow^2(S)$, 
that the intersection pairing of divisors maps only to multiples of that class:
$$\textrm{Pic}(S) \otimes \textrm{Pic}(S) 
\to \mathbb{Z} \cdot c_S \subset \chow^2 (S) ~ ,$$
and that  $c_2 (S) = 24 c_S$.

It would be interesting to see whether the fact that $\chow^2(S)$ is
generated by second Chern classes of stable vector bundles 
can be used to get a better understanding of this group.

We have shown that 
$\{c_2(M) ~|~ \textit{M is a stable vector bundle} \}$
is a set of generators for 
$\chow^2(S)$. This set is ``very big''
as we are varying arbitrarily the rank and Chern classes of the stable vector
bundles.
However we can limit this set 
even in cases where the Chow group is very big.

\begin{prop}

For every surface $S$ there is a bounded family $\mathcal{V}$ of
stable vector bundles on $S$, such that the second Chern classes of vector bundles
in $\mathcal{V}$ generate the Chow group of zero cycles in this surface.

\end{prop}

\begin{proof}
A bounded family of generating stable bundles can be constructed in various ways.
We can consider the fact that,
as $0$-cycles are formal sums of points on the surface,
then to generate the Chow group we just need to generate any (rational class of) 
single point on the surface.

We want to find a bounded family of stable vector bundles, such that their 
second Chern classes generate every point.

We can apply our construction to find a resolution of the ideal sheaf
$\ii_Z$, with $Z = \{s\}$ a point on the surface $S$. 

Following the proof of 
proposition \ref{surf}
we see that the numerical invariants chosen in the resolution
\ref{noyzm}
$$ 0 \to M_{s,V} \to V \otimes \oo_S \to \ii_s(mH) \to 0 ~,$$
\emph{i.e.} the twisting factor $m$ and the dimension of the vector space $V$ 
do not depend on the point $ p \in S$, 
but only on the ample class of the curve $C$, and  
can be fixed for all points. 

Consider, the diagonal $\Delta \subset S \times S$, and the exact sequence
$$ 0 \to \ii_{\Delta} (mH_1) \to \oo_{S\times S} (mH_1)  
\to \oo_{\Delta} (mH_1) \to 0 ~,$$
where $H_1$ is the pull-back  of $H$ to $S\times S$ through the 
projection $p_1$ of $ S\times S$ to  the first factor.

Then $I:= p_{2*} \ii_{\Delta} (mH_1)$ is a vector bundle on $S$, whose 
fiber over $s \in S$ is $H^0(S, \ii_s (mH))$.

We can consider the Grassmannian on $S$, 
$$\gr (k, I) \to S $$
where the number $k$ is $h^0 (C, \oo_C (mH))$ as in the proof of 
proposition \ref{surf}.
Then a point on $\gr (k, I) $ corresponds to a couple $(s, V)$,
where $s \in S$ and $V \subset H^0(S, \ii_s (mH))$ is a $k$-dimensional 
subsapce.

Hence we have a bounded family $\mathcal{V}= \{ M_{s,V} \}_{(s, V) \in \gr (k, I)}$.
 And 
we have shown that  second Chern classes of the stable bundles    $ M_{s,V}  $ which are in $\mathcal{V}$
(corresponding to generic $V$'s) generate the Chow group of $0$-cycles.

\end{proof}

\section{The general case}
Let now $X$ be a variety of dimension $n > 2$, with a fixed ample divisor $H$.

We want to prove the following 

\begin{thm}
\label{gen}

For every  subscheme $Z\subset X$, its ideal sheaf $\ii_Z$ admits a resolution
\begin{equation}
\label{res}
0 \to E \to P_e \to \dots \to P_1 \to P_0 \to \ii_Z \to 0
\end{equation}
where $E$ is a stable vector bundle, 
the $P_i $ are locally free sheaves of the form $ V_i \otimes \oo_X (-m_i H)$,
and $e = \dim X -2$.

\end {thm}

By passing to a multiple of $H$ we may assume that  
a generic intersection of $n-1$ sections is a smooth curve $\cou$ 
such that $g(\cou) \geqslant 1$. 
We want to prove Theorem  \ref{gen} by the same method as in the
surface case,
 \emph{i.e.} finding vector spaces 
$V_i$ that can be identified with the spaces of all global sections 
of a stable vector bundle on a smooth curve.

\subsection{Proof of the theorem}

We recall Butler's theorem for vector bundles on curves \cite{but}:

\begin{thm}[Butler]
\label{but}
Let $\cou$ be a smooth projective curve of genus $g \geqslant 1$ over
an algebraically closed field $\Bbbk$, 
and $E$ a stable vector bundle over $\cou$ with slope $\mu (E) > 2g$,
then the vector bundle 
$M_E := \mathrm{ker}( H^0 (\cou, E) \otimes \oo_{\cou}
\twoheadrightarrow E)$ 
is stable.

\end{thm}

Let us now consider a closed  sub-scheme $Z$ of codimension at least $2$.
We want to construct a sequence as in Theorem  \ref{gen}, 
which splits into short exact sequences in the following way:
$$\xymatrix@C=0.335cm@R=0.3cm{
0  \ar[r]  & E \ar[r] & P_e \ar[rr] \ar[dr] &                             & P_{e-1} \ar[r]                            
&  \dots \ar[r]   & P_2 \ar[rr] \ar[dr] &                         & P_1 \ar[rr] \ar[dr] &                   & 
P_0 \ar[r] & \ii_Z \ar[r] & 0 \\
                &             &                              & K_{e-1} \ar[ur]   &                                        
&                        &                              &     K_1 \ar[ur] &                               &K_0 \ar[ur] &
                  &                  &    
}$$
where the $K_i$ are stable sheaves on the variety $X$ which restricted to a curve $\cou$
 (an intersection of $n-1$ generic sections of $\oo_X (H)$) are stable vector bundles $M_i$, and the $P_i = V_i \otimes \oo_X (-m_i H)$ are obtained by successively lifting the space of global sections $H^0 (\cou , M_{i-1} (m_i H))$ as in the surface case. 

\begin{sloppypar}
In other words the 
$V_i \subset  H^0 (X , K_{i-1} (m_i H))$ are spaces isomorphic to $H^0 (\cou , M_i (m_i H))$ by the restriction of global sections to the curve (for the sake of clarity we should pose in the former discussion $K_{-1} :=  \ii_Z$ and $M_{-1} := \oo_{\cou}$). 
\end{sloppypar}

We remark that the stability condition is invariant 
under tensoring by a line bundle.

\begin{proof}\textbf{(Theorem \ref{gen})}

As a first step
we want to choose  $m_0$ and  \hbox{$V_0 \subset H^0(X , \ii_Z (m_0 H))$}.

Choosing $n-1$ generic\footnote{
By \emph{generic} we mean that the element 
$(s_1, \dots , s_{n-1}) \in |\oo_X(H)|^{n-1}$ is generic.} 
sections $s_1, ... , s_{n-1} \in |\oo_X(H)|$, 
gives us a filtration of $X$ by smooth sub-varieties:
$$X_0:= X \supset X_1= V(s_1) \supset X_2= V(s_1 , s_2) \supset ...
\supset X_{n-1} = \cou =  V(s_1 , \dots , s_{n-1}) .$$

Let $V \subset H^0 (X , \ii_Z (mH))$ be a subspace generating $\ii_Z (mH)$.
The restriction of  the exact sequence
$$0 \to K \to V \otimes \oo_X \to \ii_Z (mH) \to 0$$
to the hypersurface $X_1$ yields an exact sequence
$$0 \to K |_{X_1} \to V \otimes \oo_{X_1} \to \ii_Z \otimes \oo_{X_1} (mH) \to 0 ~,$$ 
due to the generality of the sections.

Restricting further we eventually obtain an exact sequence 
$$0 \to K |_{\cou} \to V \otimes \oo_{\cou} \to  \oo_{\cou} (mH) \to 0 $$
of vector bundles on the curve $\cou$.
In other words we are supposing the sequence $(s_1, \dots , s_{n-1})$ to be regular for $\ii_Z$, and such that $\cou \cap Z = \emptyset$, both of
which are open conditions. 
Furthermore, $(s_1, \dots , s_{n-1})$ being generic, we can suppose that all the 
\hbox{$\mathcal{T}or^q_{\oo_{X_i}} (\ii_Z  |_{X_{i}}, \oo_{X_{i+1}})$}
vanish, for $q>0$ and $i = 0 , \dots , n-2$:

to see this, let us fix an arbitrary locally free resolution 
$$0 \to F_s \to \dots \to F_0 \to \ii_Z \to 0$$ 
of $\ii_Z$, which splits into short exact sequences
$0 \to P_{i} \to F_{i} \to P_{i-1} \to 0$.
The sequence $(s_1, \dots , s_{n-1})$ being generic, we can suppose that it is regular for the shaves $\ii_Z , P_0 , \dots , P_{s-1}$.    Hence, from the short exact sequences above, we deduce that $\mathcal{T}or^q_{\oo_{X_i}} (\ii_Z  |_{X_{i}}, \oo_{X_{i+1}}) \cong 
\mathcal{T}or^1_{\oo_{X_i}} (P_{q-2}  |_{X_{i}}, \oo_{X_{i+1}}) = 0$.

For $m \gg 0$, we have $H^1 (X_i , \ii_Z \otimes \oo_{X_i}((m-1)H)) =0$
for every $i$. 
As in Lemma \ref{restr}, a generic 
$V \subset H^0 (X , \ii_Z (mH))$ of dimension $h^0 ( \cou , \oo_{\cou} (m))$ will map injectively to the global sections on the $X_i$:
$$\xymatrix{
& 0 \ar[d] \\ 
& H^0 (X_i , \ii_Z \otimes \oo_{X_i} ((m-1)H)) \ar[d]^{\cdot s_{i+1}} \\
V ~ \ar@{^{(}->}[r] \ar[d]_{\wr} & H^0 (X_i , \ii_Z \otimes \oo_{X_i} (mH)) \ar[d] \\
V ~ \ar@{^{(}->}[r] & H^0 (X_{i+1} , \ii_Z \otimes \oo_{X_{i+1}} (mH)) \ar[d] \\
& 0
}$$ 
until we have an isomorphism $V \tilde{\to} H^0 ( \cou , \oo_{\cou} (m))$.

So we can choose  $m_0 \gg 0$ and  $V$ generating
$\ii_Z (m_0 H)$ such that  the kernel $K_0$ of 
\hbox{$V \otimes \oo_X (-m_0) \to \ii_Z$} 
is stable (since it's stable on the curve $\cou$ which is a complete intersection of $n-1$ sections of $H$),
but $K_0$ is, in general, not locally free.

As we have chosen $(s_1, \dots , s_{n-1})$ such that  
$\mathcal{T}or^q_{\oo_{X_i}} (\ii_Z  |_{X_{i}}, \oo_{X_{i+1}}) = 0$
for $q>0$ and $i = 0 , \dots , n-2$, we deduce from the sequence
$$0 \to K_0 \to V \otimes \oo_X (-m_0) \to \ii_Z \to 0$$ that also the 
$\mathcal{T}or^q_{\oo_{X_i}} (K_0 |_{X_i} , \oo_{X_{i+1}})$
vanish, for $q>0$ and $i = 0 , \dots , n-2$. In particular, the sequence 
$(s_1, \dots , s_{n-1})$ is $K_0$-regular.

Repeating the argument, we  obtain, tensoring  $K_0$ by $H$ enough times, exact sequences:
$$o \to K_1(m_1 H) |_{X_i} \to V_1 \otimes \oo_{X_i}  \to K_0 (m_1 H) |_{X_i} \to 0 ~.$$

Again, we can suppose that $H^1 (X_i , K_0 \otimes \oo_{X_i} (m_1 H) ) = 0$ 
and lift the vector space $H^0 (\cou , K_0(m_1 H) |_{\cou})$ 
on a generic generating space 
$V_1 \subset H^0 (X , K_0(m_1 H))$.
Butler's theorem tells us that the vector bundle ${K_1} |_{\cou}$, satisfying 
$$ 0 \to K_1(m_1 H) |_{ \cou} \to H^0(\cou , K_0(m_1 H) |_{\cou}) \otimes \oo_{\cou}  
\to K_0 (m_1 H) |_{\cou} \to 0 ~,$$
is a stable vector bundle (for $m_1 \gg 0$), because ${K_0} |_{\cou}$ is stable and locally free.

So we can continue and find the resolution (\ref{res}), where we remark that if 
\hbox{$e \geqslant n-2$},
$E$ is a vector bundle because $X$ is smooth and so has cohomological  dimension 
$n = \dim X$, and it is stable because it is so on the curve $\cou$.

\end{proof}

\subsection{Stable vector bundles as generators}

We can apply then this result to calculate the Chern class and character of $\ii_Z$;
we know that in general for any sheaf $\fais$ and any resolution 
$P^{\bullet} \to \fais$ by vector bundles, its Chern  character is 
$ ch(\fais) = \sum {(-1)^i} ch( P^i)$.

\begin{cor}

A set of generators of $\chowqx$, as a group, is 
$$\{ ch(E) | E \textrm{ stable vector bundle}\} ~.$$

\end{cor}

\vspace{0cm}

\begin{proof}
From the resolution (\ref{res}) we have:
$$ ch(\ii_Z) = (-1)^{e +1} ch(E) + 
\sum_{i = 0}^{e} (-1)^i \dim V_i \cdot  ch(\oo_X (- m_i H)) ~.$$

From the theorem of Grothendieck-Riemann-Roch (cf. \cite{grot}) we know that 
$$ch(\ii_Z) = 1 - ch( \oo_Z)= 1 - [Z] + \textrm{higher order terms}$$
so applying our result to the higher order terms, 
we see that we can express $[Z]$ as 
a sum  of Chern characters of stable vector bundles. 

\end{proof}

In order to have the same results in the $K$-theory and the derived category we will use the following 

\begin{lem}

Any coherent sheaf $\fais$on $X$ admits a filtration 
$0 = \fais_0 \subset \fais_1 \subset \dots \subset \fais_{\ell} = \fais$ where each quotient 
$ \fais_i / \fais_{i-1}$  admits a polystable 
resolution.

\end{lem}

\begin{proof}
Consider at first a torsion sheaf $\tais$:
it has  then a filtration 
$0 =  \tais_0 \subset \tais_1 \subset \dots \subset \tais_{\ell} = \tais$, 
where every quotient $ \tais_i  / \tais_{i-1}$ is of the form $\oo_{Z_i} (-m H)$, for cycles $Z_i~$. 
Hence $\tais$  admits such a filtration.

A torsion free sheaf $\fais$ admits an extension 
$$0 \to V \otimes \oo_X (-m) \to \fais \to \frac{\fais} {V \otimes \oo_X (-m)} \to 0 ~,$$
where $m\gg0$,  $V \subseteq H^0 ( X , \fais (m))$ is  the  subspace generated by 
$R$ generically independent sections of $\fais (m)$, $R$ is the generic rank of $\fais$,
and ${\fais} / (V \otimes \oo_X (-m))$ is a torsion sheaf. 
Hence taking the pull-back to $\fais$ of the torsion sheaf filtration, we get the requested filtration.

Finally, any coherent sheaf fits into an extension with its torsion and torsion free parts:
$$ 0 \to \tais (\fais) \to \fais \to \fais / \tais (\fais) \to 0 ~,$$
so we can take the filtration for $\tais (\fais)$ and the pull-back to $\fais$ of the filtration for 
$\fais / \tais (\fais)$.

\end{proof}

The following result is an immediate consequence:

\begin{cor}

The Grothendieck ring $K(X)$ is generated, as a group, by the classes of stable vector bundles.

\end{cor}

\begin{rmk} 
\label{notfree}
Every torsion free sheaf admits a (unique) Harder-Narashiman filtration, whose
quotients are semistable sheaves (not necessarely locally free).
And every semistable sheaf admits a (non unique) filtration with
stable quotients. Mixing those two kinds of filtrations 
we obtain a filtration with stable quotients of any torsion free sheaf.

Hence,
it can be easily proven that the  class in $K(X)$ of any coherent sheaf 
$\fais$
is obtained as a sum of classes 
of stable not necessarily locally free sheaves.
In fact we can
construct an exact sequence
$0 \to K \to V \otimes \oo_X(-mH) \to \fais \to 0$, 
and take the filtration of the torsion free sheaf $K$, whose quotients
 are stable not necessarily locally free sheaves. (The same argument
 holds for the Chow group).
\end{rmk}

For what concerns the derived category, let $\dboun{X}$ be the bounded derived category of coherent sheaves on X.
We will identify, as usual, any coherent sheaf $\fais$ to the object $(0 \to \fais \to 0) \in \dboun{X}$
 concentrated in degree 0.

\begin{dfn}
We say that a triangulated subcategory 
$\mathcal{D} \subseteq \dboun{X}$, is generated by a
family of objects $\mathcal{E} \subseteq \dboun{X}$, if it is 
 the smallest triangulated full subcategory of $\dboun{X}$, stable under isomorphisms, 
which contains 
$\mathcal{E}$. We will denote it by $<\mathcal{E}>$.

\end{dfn}

It is easy to prove the following lemmas:

\begin{lem}

Let $\mathcal{E}$ be a family of objects of $\dboun{X}$. 
If $<\mathcal{E}>$   contains two coherent sheaves $\fais_1$ and $\fais_2$, then
it contains all their  extensions.

\end{lem}

\begin{lem}

Let $\mathcal{E}$ be a family of objects of $\dboun{X}$.
If $<\mathcal{E}>$ contains every coherent sheaf, 
then $<\mathcal{E}> = \dboun{X}$.

\end{lem}

As in the case of the Grothendieck group, we get immediately the following 

\begin{cor}

The bounded derived category $\dboun{X}$ is generated by the family of stable vector bundles.

\end{cor}

\bibliographystyle{amstronzo}
\bibliography{chow}

 \vspace{1cm}

 \begin{flushright}

 \small{

 Ernesto Carlo \textsc{Mistretta}

 \texttt{ernesto@math.jussieu.fr}

 \textsc{Institut de Mathématiques de Jussieu}

 Équipe de Topologie et Géométrie Algébriques

 175, rue de Chevaleret

 75013 Paris}

 \end{flushright}

\end{document}